\documentclass{amcjou}

\usepackage{amssymb}

\usepackage{graphicx}

\usepackage{hyperref}
\usepackage{cleveref}

\makeatletter
\newcommand{\noprelistbreak}{\@nobreaktrue\nopagebreak\medskip} 
\makeatother

\newcommand{\pref}[1]{(\ref{#1})}
\newcommand{\csee}[1]{(see \cref{#1})}

\newcommand{\fullcref}[2]{\cref{#1}\pref{#1-#2}}
\newcommand{\fullcsee}[2]{(see \fullcref{#1}{#2})}

\newcommand{\integer}{\mathbb{Z}}
\newcommand{\quot}{\overline}

\DeclareMathOperator{\Cay}{Cay}

\newcommand{\connsum}{\mathbin{\#}}
\newcommand{\voltage}{\Pi}
\newcommand{\Voltages}[1]{\mathcal{V}_{#1}}

\newcommand{\halfit}[1]{%
	\mathchoice
		{\setbox0\hbox{\hskip1pt$\displaystyle#1$\hskip1pt}\copy0%
		\kern-\wd0%
		\hbox{\vrule height 2.5pt depth-2pt width \wd0}}%
		{\setbox0\hbox{\hskip1pt$\textstyle#1$\hskip1pt}\copy0%
		\kern-\wd0%
		\hbox{\vrule height 2.5pt depth-2pt width \wd0}}%
		{\setbox0\hbox{\hskip0.5pt$\scriptstyle#1$\hskip0.5pt}\copy0
		\kern-\wd0%
		\hbox{\vrule height 1.5pt depth-1pt width \wd0}}%
		{\setbox0\hbox{\hskip0.5pt$#1$\hskip0.5pt}\copy0%
		\kern-\wd0%
		\hbox{\vrule height 1.5pt depth-1pt width \wd0}}%
	  }

\newcommand{\MR}[1]{MR\,#1}

\numberwithin{equation}{section}

\makeatletter 
\swapnumbers
\def\swappedhead#1#2#3{%
  {\normalfont(\thmnumber{#2})}
  \thmname{\@ifnotempty{#2}{~}#1}%
  \thmnote{ {\the\thm@notefont#3}}} 
\def\thmhead@plain#1#2#3{%
  \thmname{#1}\thmnumber{\@ifnotempty{#1}{ }\@upn{#2}}%
  \thmnote{ {\the\thm@notefont#3}}} 
\let\thmhead\thmhead@plain
\makeatother

\newtheorem{prop}[equation]{Proposition}
\newtheorem{thm}[equation]{Theorem}
\newtheorem{lem}[equation]{Lemma}
\newtheorem{cor}[equation]{Corollary}

\newtheorem{FGL}[equation]{Lemma}

\theoremstyle{definition}
\newtheorem{rem}[equation]{Remark}
\newtheorem{notation}[equation]{Notation}
\newtheorem{defn}[equation]{Definition}
\newtheorem*{ack}{Acknowledgments}

\crefformat{FGL}{the Factor Group Lemma~#2(#1)#3}

 \newcounter{case}
 \newenvironment{case}[1][\unskip]{\refstepcounter{case}\bf
 \medskip \noindent Case \thecase\ #1. \it}{\unskip\upshape}
 \renewcommand{\thecase}{\arabic{case}}
\crefformat{case}{Case~#2#1#3}

 \newcounter{subcase}
 
\crefformat{subcase}{Subcase~#2#1#3}
\renewcommand{\thesubcase}{\thecase.\arabic{subcase}}
\numberwithin{subcase}{case}

\crefformat{figure}{Figure~#2#1#3}
\crefmultiformat{figure}{Figures~#2#1#3}{ and~#2#1#3}{, #2#1#3}{, and~#2#1#3}

\crefformat{section}{Section~#2#1#3}

\makeatletter
\renewenvironment{frontmatter}
{\thispagestyle{amctitle}%
\setcounter{page}{\@startpage}%
\vskip 10pt%
\centering}
{\vskip 20pt%
\blfootnote{\raggedright \ifnum\@authorcount=1\textit{E-mail address:}\else\textit{E-mail addresses:}\fi ~\@emails}%
}
\renewcommand{\ps@amctitle}{%
  \renewcommand\@oddhead{}
  \let\@evenhead\@oddhead
  \renewcommand\@evenfoot{\hfil}
  \let\@oddfoot\@evenfoot
}
\def\@oddrunninghead{Cayley graphs on nilpotent groups with cyclic commutator subgroup are hamiltonian}
\def\@evenrunninghead{Ebrahim Ghaderpour and Dave Witte Morris}
\makeatother

\begin{document}

\begin{frontmatter}

\titledata{Cayley graphs on nilpotent groups with cyclic commutator subgroup are hamiltonian}{}
{(version of 26 November 2011)} 

\authordata{Ebrahim Ghaderpour}
{Department of Mathematics and Computer Science,
University of Lethbridge, Lethbridge, Alberta, T1K~3M4, Canada}
{Ebrahim.Ghaderpoor@uleth.ca}
{}

\authordata{Dave Witte Morris}
{Department of Mathematics and Computer Science,
University of Lethbridge, Lethbridge, Alberta, T1K~3M4, Canada}
{Dave.Morris@uleth.ca, http://people.uleth.ca/\!$\sim$dave.morris/}
{}


\begin{abstract}
We show that if $G$ is any nilpotent, finite group, and the commutator subgroup of~$G$ is cyclic, then every connected Cayley graph on~$G$ has a hamiltonian cycle.
\end{abstract}

\end{frontmatter}

\section{Introduction}

It has been conjectured that every connected Cayley graph has a hamiltonian cycle. See \cite{CurranGallian-survey, PakRadoicic-survey, Witte-CayDiags, WitteGallian-survey} for  references to some of the numerous results on this problem that have been proved in the past forty years, including the following theorem that is the culmination of papers by Maru\v si\v c \cite{Marusic-HamCircCay}, Durnberger \cite{Durnberger-semiprod,Durnberger-prime}, and  Keating-Witte \cite{KeatingWitte}:

\begin{thm}[(D.\,Maru\v si\v c, E.\,Durnberger, K.\,Keating, and D.\,Witte, 1985)] \label{KeatingWitteThm}
Let $G$ be a nontrivial, finite group. If the commutator subgroup $[G,G]$ of~$G$ is cyclic of prime-power order, then every connected Cayley graph on~$G$ has a hamiltonian cycle.
\end{thm}

It is natural to try to prove a generalization that only assumes the commutator subgroup is cyclic, without making any restriction on its order, but that seems to be an extremely difficult problem: at present, it is not even known that all connected Cayley graphs on dihedral groups are hamiltonian. (See \cite{AlspachChenDean, AlspachZhang} and \cite[Cor.~5.2]{Witte-CayDiags} for the main results that have been proved for dihedral groups.) In this paper, we replace the assumption on the order of $[G,G]$ with the rather strong assumption that $G$ is nilpotent:

\begin{thm} \label{NilpG'Cyclic}
Let $G$ be a nontrivial, finite group. If $G$ is nilpotent, and the commutator subgroup of~$G$ is cyclic, then every connected Cayley graph on~$G$ has a hamiltonian cycle.
\end{thm}

The proof of this \namecref{NilpG'Cyclic} is based on a variant of the method of D.\,Maru\v si\v c \cite{Marusic-HamCircCay} that established \cref{KeatingWitteThm} (cf.\ \cite[Lem.~3.1]{KeatingWitte}). 

\begin{rem} \label{Previous}
Here are some previous results on the hamiltonicity of the Cayley graph $\Cay(G;S)$ when $G$ is nilpotent:
\noprelistbreak
\begin{enumerate}

\item Assume $G$ is nilpotent, the commutator subgroup of~$G$ is cyclic, and $\#S = 2$. Then a hamiltonian cycle in $\Cay(G;S)$ was found in \cite[\S6]{KeatingWitte} \csee{2genNilpG'Cyclic}. The present paper generalizes this by eliminating the restriction on the cardinality of the generating set~$S$.

\item For Cayley graphs on nilpotent groups (without any assumption on the commutator subgroup), it was recently shown that if the valence is at most~$4$, then there is a hamiltonian \emph{path} (see \cite{Morris-2GenNilp}). 

\item Every nilpotent group is a direct product of $p$-groups. For $p$-groups, it is known that every Cayley graph has a hamiltonian cycle (\cite{Witte-pgrp}, see \cref{3grp}).  Unfortunately, we do not know how to extend this to direct products.

\item  \label{Previous-abel}
Every abelian group is nilpotent. It is well known (and easy to prove) that Cayley graphs on abelian groups always have hamiltonian cycles. In fact, they are usually hamiltonian connected (see \cite{ChenQuimpo-hamconn}).

\end{enumerate}
\end{rem}

\begin{ack}
We thank Dragan Maru\v si\v c and Mohammad Reza Salarian for their comments that encouraged this line of research.
\end{ack}

\section{Assumptions, notation, and outline of the proof} \label{AssumpsSect}

We begin with some standard notation:

\begin{notation} 
Let $G$ be a group, and let $S$ be a subset of~$G$.
	\noprelistbreak
	\begin{itemize}
	\item $\Cay(G;S)$ denotes the \emph{Cayley graph} of~$G$ with respect to~$S$. Its vertices are the elements of~$G$, and there is an edge joining~$g$ to~$gs$ for every $g \in G$ and $s \in S$.
	\item $G' = [G,G]$ denotes the commutator subgroup of~$G$.
	\item $S^r = \{\, s^r \mid s \in S\,\}$ for any $r \in \integer$. 
	\item $S^{\pm1} = S \cup S^{-1}$. 
	\end{itemize}
Note that if $S$ happens to be a cyclic subgroup of~$G$, then $S^r$ is a subgroup of~$S$.
\end{notation}

We now fix notation designed specifically for our proof of \cref{NilpG'Cyclic}:

\begin{notation} \ 
	\noprelistbreak
	\begin{itemize}
	\item $G$ is a nilpotent, finite group,
	\item $N$ is a cyclic, normal subgroup of $G$ that contains~$G'$,
	\item $g \mapsto \quot{g}$ is the natural homomorphism from~$G$ to $G/N = \quot{G}$,
	\item $S = \{\sigma_1,\sigma_2,\ldots,\sigma_\ell\}$ is a subset of~$G$, such that 
		\begin{itemize}
		\item $\quot{S}$ is a minimal generating set for $\quot{G}$, 
		and
		\item $\ell = \#S = \#\quot{S} \ge 2$,
		\end{itemize}
	\item $S_k = \{\, \sigma_i \mid i \le k \,\}$ for $1 \le k \le \ell$,
	\item $G_k = \langle S_k \rangle N$,
	and
	\item $m_k = | \quot{G_k} : \quot{G_{k-1}} |$.

	\end{itemize}
\end{notation}

 \begin{defn} \ 
	\noprelistbreak
	\begin{itemize}

	\item If $(s_i)_{i=1}^n$ is a sequence of elements of~$S^{\pm1}$, and $\quot{g} \in \quot{G}$, we use $\quot{g}(s_i)_{i=1}^n$ to denote the walk in $\Cay(\quot{G};\quot{S})$ that visits (in order) the vertices
	$$ \quot{g}, \quot{g s_1}, \quot{g s_1 s_2}, \ldots, , \quot{g s_1 s_2 \cdots s_n} .$$

	\item  If $C = \quot{g}(s_i)_{i=1}^n$ is any oriented cycle in $\Cay(\quot G; \quot S)$, its \emph{voltage} is $\prod_{i=1}^n s_i$. This is an element of~$N$, and it may be denoted $\voltage C$.
	\item For $S_0 \subset S$, we say the walk $\quot{g}(s_i)_{i=1}^n$ \emph{covers $S_0^{\pm1}$} if it contains an oriented edge labeled~$s$ and a {\upshape(}different\/{\upshape)} oriented edge labeled~$s^{-1}$, for every $s \in S_0$. (That is, there exist $i,j$ with $i \neq j$, such that $s_i = s$ and $s_j = s^{-1}$.)
	\item $\Voltages{k}$ is the set of voltages of oriented hamiltonian cycles in $\Cay(\quot{G_k}; \quot{S_k})$ that cover~$S_k^{\pm1}$.
	\end{itemize}
 \end{defn}
 
The following well-known, elementary observation is the foundation of our proof:

\begin{FGL}[(``Factor Group Lemma'' {\cite[\S2.2]{WitteGallian-survey}})] \label{FGL}
Suppose
 \begin{itemize}
 \item $N$ is a cyclic, normal subgroup of~$G$,
 \item $C = \quot{g}(s_i)_{i=1}^n$ is a hamiltonian cycle in $\Cay(\quot{G};\quot{S})$,
 and
 \item the voltage $\voltage C$ generates~$N$.
 \end{itemize}
 Then $(s_1,\ldots,s_n)^{|N|}$ is a hamiltonian cycle in $\Cay(G;S)$.
 \end{FGL}

With this in mind, we let $N = G'$, and we would like to find a hamiltonian cycle in $\Cay(\quot{G}; \quot{S})$ whose voltage generates~$N$. In almost all cases, we will do this by induction on $\ell = \#S$, after substantially strengthening the induction hypothesis. Namely, we consider the following assertion ($\alpha_k^\epsilon$) for $2 \le k \le \ell$ and $\epsilon \in \{1,2\}$:  
	{\renewcommand{\theequation}{$\alpha_k^\epsilon$}
	\begin{align}
	 \label{AlphaDefn}
	\begin{matrix}
	\text{there exists $h_k \in N$, such that, for every $x \in N$,}  \\
	\text{$\bigl( \Voltages{k} \cap  (G_k')^\epsilon h_k \bigr) x$ contains a generator of} \\
	\text{a subgroup of~$N$ that contains~$(G_k')^\epsilon$.} 
	\end{matrix}
	\end{align}
	}%
	\addtocounter{equation}{-1}%
For $\epsilon = 2$, we also consider the following slightly stronger condition, which we call $\alpha_k^{2+}$:
	{\renewcommand{\theequation}{$\alpha_k^{2+}$}
	\begin{align}
	 \label{Alpha+Defn}
	\text{$\alpha_k^2$ holds, and $\langle h_k, (G_k')^2\rangle$ contains $G_k'$.} 
	\end{align}
	}%
	\addtocounter{equation}{-1}%

\begin{lem} \label{Alpha->Spiralling}
Let $N = G'$. If either $\alpha_\ell^1$ or $\alpha_\ell^{2+}$ holds, then there is a hamiltonian cycle in $\Cay(\quot{G}; \quot{S})$ whose voltage generates~$N$.
\end{lem}

\begin{proof}
Note that $G_\ell' = G' = N$. Since $\Voltages{\ell}$ consists of voltages of hamiltonian cycles in $\Cay(\quot{G}; \quot{S})$, it suffices to find an element of $\Voltages{\ell}$ that generates~$G_\ell'$.

If we assume $\alpha_\ell^1$, then the desired conclusion is immediate, by taking $x = e$ in that assertion.

Similarly,  if we assume $\alpha_\ell^{2+}$, then taking $x = e$ in~$\alpha_\ell^2$  tells us that some element~$\gamma$ of~$\Voltages{\ell} \cap (G_\ell')^2 h_\ell$ generates a subgroup that contains~$(G_\ell')^2$. Then, since $\gamma \in (G_\ell')^2 h_\ell$, and $\langle h_\ell, (G_\ell')^2\rangle$ contains $G_\ell'$, we have
	\begin{align*}
	\langle \gamma \rangle &= \langle \gamma, (G_\ell')^2\rangle = \langle h_\ell, (G_\ell')^2\rangle \supset G_\ell' = N
	.  \qedhere \end{align*}
\end{proof}

\begin{rem} \ \label{Alpha<>Alpha}
	\noprelistbreak
	\begin{enumerate}
	\item \label{Alpha<>Alpha-odd}
	If $|G_k'|$ is odd, then $(G_k')^2 = G_k'$, so we have $\alpha_k^1 \Leftrightarrow \alpha_k^2 \Leftrightarrow \alpha_k^{2+}$ in this case.  Thus, the parameter~$\epsilon$ is only of interest when $|G'|$ is even.

	\item It is not difficult to see that $\alpha_k^1 \Rightarrow \alpha_k^{2+}$, but we do not need this fact.
	\end{enumerate}
\end{rem}

Our proof of $\alpha_\ell^1$ or $\alpha_\ell^{2+}$ is by induction on~$k$. Here is the outline:
	\noprelistbreak
	\begin{enumerate} \renewcommand{\theenumi}{\Roman{enumi}}
	\item We prove a base case of the induction:
	 $\alpha_2^2$ is usually true \csee{KW6-2}.
	 \item We prove an induction step:
	 under certain conditions,
	 $\alpha_k^1 \Rightarrow \alpha_{k+1}^1$ and $\alpha_k^{2+} \Rightarrow \alpha_{k+1}^{2+}$  \csee{Alpha->Alpha}.
	 \item We prove $\alpha_\ell^1$ or $\alpha_\ell^{2+}$ is usually true, by bridging the gap between $\alpha_2^2$ and either $\alpha_3^1$ or $\alpha_3^{2+}$, and then applying the induction step \csee{SpiralG'Odd,SpiralG'Even}.
	\end{enumerate}

Here is a detailed explanation of how our results combine with the main result of \cite{Witte-pgrp} to establish the main theorem:

\begin{proof}[\normalfont\textbf{Proof of \cref{NilpG'Cyclic}}]
We may assume:
	\begin{itemize}
	\item $\ell \ge 3$, for otherwise \cref{2genNilpG'Cyclic} applies.
	\item $|G|$ is not a power of~$3$, for otherwise \cref{3grp}  applies.  
	\end{itemize}

Let $S$ be a minimal generating set of~$G$, and let $N = G'$. Note that $\quot{S}$ is a minimal generating set of~$\quot{G}$ (because $G'$ is contained in the Frattini subgroup $\Phi(G)$, cf.\ \cite[Cor.~10.3.3]{Hall-ThyOfGrps}).

We claim there is a hamiltonian cycle in $\Cay(\quot{G};\quot{S})$ whose voltage generates~$G'$. While proving this, there is no harm in assuming that $|G'|$ is square-free \csee{FreeLunch}. Also note that, since $|G|$ is not a power of~$3$, we cannot have $|G'| = |\quot{s}| = 3$ for all $s \in S$. Then, by applying either \cref{SpiralG'Odd} or \cref{SpiralG'Even} (depending on the parity of $|G'|$), we obtain either $\alpha_\ell^1$ or~$\alpha_\ell^{2+}$. Each of these yields the desired hamiltonian cycle in $\Cay(\quot{G};\quot{S})$ \csee{Alpha->Spiralling}. 

Now that the claim has been verified, \cref{FGL} provides a hamiltonian cycle in $\Cay(G;S)$.
\end{proof}

\section{Preliminaries}

\subsection{Results from \cite{KeatingWitte} and \cite{Witte-pgrp}}

The following result from \cite{Witte-pgrp} allows us to assume $G$ is not a $3$-group. (Since we always assume that $G'$ is cyclic, a short proof of the special case we need can be found in \cite[Thm.~6.1]{Witte-CayDiags}.)

\begin{prop}[(Witte \cite{Witte-pgrp})] \label{3grp}
If $|G|$ is a power of some prime~$p$, then every connected Cayley graph on~$G$ has a hamiltonian cycle.
\end{prop}

The following simple observation usually allows us to assume $|N|$ is square-free.

\begin{lem}[{}{\cite[Lem.~3.2]{KeatingWitte}}] \label{FreeLunch}
Let $\underline{G} = G/\Phi(N)$, where $\Phi(N)$ is the Frattini subgroup of~$N$ \cite[\S10.4]{Hall-ThyOfGrps}. Then:
	\noprelistbreak
	\begin{enumerate}
	\item $|\underline{N}|$ is square-free,
	and
	\item if there is a hamiltonian cycle in $\Cay( \underline{G}/\underline{N} ; S)$ whose voltage generates~$\underline{N}$, then there is a hamiltonian cycle in $\Cay( G/N ; S)$ whose voltage generates~$N$.
	\end{enumerate}
\end{lem}

\begin{lem}[(Keating-Witte {\cite[Case~6.1]{KeatingWitte}})] \label{TwoGen}
If $|\quot{G_2}|$ is even, then $\Cay \bigl( \quot{G_2}; \quot{S_2} \bigr)$ has a hamiltonian cycle whose voltage is a generator of~$G_2'$.
\end{lem}

\begin{proof}
For the reader's convenience, we reproduce the gist of the argument, since it is very short. We may assume $|\sigma_1|$ is even (by interchanging $\sigma_1$ and~$\sigma_2$ if necessary). For convenience, let $n = |\sigma_1|$ and $m = m_2$. Then 
	$$ \bigl( \sigma_2^{m-1},(a,\sigma_2^{-(m-2)},a, \sigma_2^{m-2})^{(n-2)/2}, a, \sigma_2^{-(m-1)}, \sigma_1^{-(n-1)} \bigr) $$
is a hamiltonian cycle in $\Cay( \quot{G};\quot{S})$ whose voltage is $[\sigma_1,\sigma_2]$, which generates $G_2'$ \csee{CommsGen}.
\end{proof}

The following result allows us to assume $\ell \ge 3$.

\begin{prop}[(Keating-Witte {\cite[\S6]{KeatingWitte}})] \label{2genNilpG'Cyclic}
If $\ell = 2$ and $N = G'$, then $\Cay \bigl( G; S \bigr)$ has a hamiltonian cycle.
\end{prop}

\begin{proof}
For the reader's convenience, we point out how to derive this from results proved in this paper (and \cref{3grp}).

We may assume $|G/G'|$ is odd, for otherwise a hamiltonian cycle is obtained by combining \cref{TwoGen} with \cref{FGL}.
We may also assume that $|G|$ is not a power of~$3$, for otherwise \cref{3grp} applies. This implies it is not the case that $|\quot{s}| = 3$ for every $s \in S$.

If $|G'|$ is square-free, then \cref{KW6-2} tells us that $\alpha_2^2$ is true. Since $|G'|$ is odd, this implies that $\alpha_2^1$ is true \fullcsee{Alpha<>Alpha}{odd}. So \cref{FGL} provides a hamiltonian cycle in $\Cay(G;S)$ \csee{Alpha->Spiralling}, and \cref{FreeLunch} tells us there is a hamiltonian cycle even without the assumption that $|G'|$ is square-free.
\end{proof}

\subsection{Remarks on voltage}

 \begin{rem}
 By definition, it is clear that all translates of~$C$ have the same voltage.
 That is,
 	$$ \voltage \bigl( \quot{g}(s_i)_{i=1}^n \bigr) = \voltage \bigl( (s_i)_{i=1}^n \bigr) .$$
\end{rem}

 \begin{rem}
If $|N|$ is square-free (which is usually the case in this paper), then $N$ is contained in the center of~$G$ (because $|N|$ is the direct product of normal subgroups of prime order, and it is well known that those are all in the center \cite[Thm.~4.3.4]{Hall-ThyOfGrps}). In this situation, the voltage of a cycle is independent of the starting point that is chosen for its representation. That is, if $(t_i)_{i=1}^n$ is a cyclic rotation of~$(s_i)_{i=1}^n$, so there is some $r \in \{0,1,2,\ldots,n\}$ with $t_i = s_{i+r}$ for all~$i$ (where subscripts are read modulo~$n$), then
 	$$
	 \voltage  (t_i)_{i=1}^n 
	 = s_{r+1} s_{r+2} \cdots s_n \, s_1 s_2 \cdots s_r
	= (s_1s_2\cdots s_r)^{-1} \bigl( \voltage (s_i)_{i=1}^n \bigr) s_1s_2\cdots s_r
	= \voltage (s_i)_{i=1}^n
	,$$
because $\voltage (s_i)_{i=1}^n  \in N \subset Z(G)$.
 \end{rem}

\subsection{Elementary facts about cyclic groups of square-free order}

When we want to show that some subgroup~$H$ of~$N$ contains some other subgroup~$K$, the  following observation often allows us to assume $K = N$ (by modding out~$K^\perp$), which means we wish to prove $H = N$. 

\begin{lem} \label{ContainsProj}
Assume $|N|$ is square-free, and $H$ and~$K$ are two subgroups of~$N$.
Then:
	\noprelistbreak
	\begin{enumerate}
	\item \label{ContainsProj-Kperp}
	There is a unique subgroup~$K^\perp$ of~$N$, such that $N = K \times K^\perp$.
	\item \label{ContainsProj-normal}
	$K^\perp$ is a normal subgroup of~$G$.
	\item \label{ContainsProj-subset}
	$K \subseteq H$ iff $\underline{H} = \underline{N}$ in $\underline{G} = G/K^\perp$.
	\end{enumerate}
\end{lem}

\begin{proof}
\pref{ContainsProj-Kperp}~Since $N$ is cyclic, it has a unique subgroup of any order dividing~$|N|$; let $K^\perp$ be the subgroup of order $|N/K|$. Since $|N|$ is square-free, we have $\gcd \bigl( |K|, |K^\perp| \bigr) = 1$, so $N = K \times K^\perp$.

\pref{ContainsProj-normal}~It is well known that every subgroup of a cyclic, normal subgroup is normal (because no other subgroup of~$N$ has the same order).

\pref{ContainsProj-subset}~We prove only the nontrivial direction. Since $\underline{H} = \underline{N}$, we know that $|K| = |\underline{N}|$ is a divisor of~$|H|$. So $|H|$ has a subgroup whose order is~$|K|$. Since $K$ is the only subgroup of~$N$ with this order, we must have $K \subseteq H$. 
\end{proof}

\begin{lem} \label{GCD(consecutive)}
Suppose 
	\noprelistbreak
	\begin{itemize}
	\item $\gamma$ is a generator of~$N$,
	\item $x \in N$,
	and
	\item $a \ge \max \bigl( |N|, 5 \bigr)$.
	\end{itemize}
Then, for some $i$ with $1 \le i \le \lfloor (a-1)/2 \rfloor$, we have $N^2 \subseteq \langle \gamma^{-2i} x \rangle$.
\end{lem}

\begin{proof}
Write $x=\gamma^h$, where $1\leq h \leq |N|$, choose $r \in \{1,2\}$ such that $h - r$ is even, and let 
	$$i = \begin{cases}
	\hfil r & \text{if $h \in \{1,2\}$}, \\
	(h-r)/2 & \text{if $h > 2$} 
	. \end{cases}$$
Then $h - 2i \in \{\pm r\} \subset \{\pm1, \pm2\}$, so $N^2 \subseteq\langle\gamma^{h-2i} \rangle= \langle \gamma^{-2i} x \rangle$.
\end{proof}

\begin{lem} \label{ContainsN2}
If 
	\noprelistbreak
	\begin{itemize}
	\item $N$ is a cyclic group of square-free order,
	\item $m \ge |N|$,
	\item $k \ge 2$,
	\item $T = \{\gamma_1,\ldots,\gamma_k\}$ generates~$N$, 
	and
	\item $h \in N$,
	\end{itemize}
then we may choose a sequence $(j_i)_{i=1}^{m-1}$ of elements of $\{1,2,\ldots,k\}$, and $\gamma_i^* \in \{\gamma_{j_i}^{\pm1}\}$ for each~$i$, such that  $\gamma_{i+1}^* = \gamma_i^*$ whenever $j_{i+1} = j_i$, and
	\begin{align} \label{ContainsN2-display}
	 \text{$\langle h \gamma_1^* \gamma_2^* \cdots \gamma_{m-1}^* \rangle$ contains $N^2$} 
	 . \end{align}
Furthermore, if either 
	\begin{enumerate}
	\item \label{ContainsN2-odd}$|N|$ is odd, 
	or 
	\item  \label{ContainsN2-notsame}
	the elements of~$T$ are not all in the same coset of $N^2$,
	\end{enumerate}
then $\gamma_1^*,\ldots,\gamma_{m-1}^*$ can be chosen so that \pref{ContainsN2-display} holds with $N$ in the place of~$N^2$.
\end{lem}

\begin{proof} 
%
Let us assume $|N| > 3$. (The smaller cases are very easy to address individually.)

We begin by finding $\gamma_1^*, \gamma_2^*, \ldots, \gamma_{m-1}^* \in T^{\pm1}$, such that $\langle h 
\gamma_1^* \gamma_2^* \cdots \gamma_{m-1}^* \rangle$ contains~$N^2$ (or $N$, if appropriate), but without worrying about the requirement that $\gamma_{i+1}^* = \gamma_i^*$ whenever $j_{i+1} = j_i$. 

Assume, for the moment, that the Cayley graph $\Cay(N;T)$ is not bipartite. (In other words, assume that either \pref{ContainsN2-odd} or~\pref{ContainsN2-notsame} holds.) Also, let $\gamma$ be a generator of~$N$, and assume $h^{-1}\gamma \neq e$ (by replacing $\gamma$ with its inverse, if necessary). Then, since $\Cay(N;T)$ is not bipartite, there is a walk $(\gamma_i^*)_{i=1}^r$ from~$e$ to~$h^{-1}\gamma$, such that $r \equiv m-1 \pmod{2}$. With a bit of care, we can also ensure that $r < |N|$, so $r \le m-1$. Then
	$$ 
	\text{$h \gamma_1^* \gamma_2^* \cdots \gamma_r^* (\gamma_1 \gamma_1^{-1})^{(m-1-r)/2} = \gamma$ generates~$N$,} $$
as desired.

Now suppose  $\Cay(N;T)$ is bipartite. Let $(N^2)^\perp$ be the subgroup of order~$2$ in~$N$, and let $\underline{N} = N/(N^2)^\perp$. Then $|\underline{N}|$ is odd, so $\Cay( \underline{N}; T)$ is certainly not bipartite. Therefore, the preceding paragraph provides $\gamma_1^*, \gamma_2^*, \ldots, \gamma_{m-1}^* \in T^{\pm1}$, such that $\langle h 
\gamma_1^* \gamma_2^* \cdots \gamma_{m-1}^* \rangle = \underline{N}$. This implies that  $\langle h \gamma_1^* \gamma_2^* \cdots \gamma_{m-1}^* \rangle$ contains~$N^2$ \csee{ContainsProj}.

\medbreak

To complete the proof, we modify the above sequence $\gamma_1^*, \gamma_2^*, \ldots, \gamma_{m-1}^*$ to satisfy the condition that  $\gamma_{i+1}^* = \gamma_i^*$ whenever $j_{i+1} = j_i$. First of all, since $N$ is commutative, we may collect like terms, and thereby write 
	$$\gamma_1^* \gamma_2^* \cdots \gamma_{m-1}^* = \gamma_1^{m_1} \gamma_2^{m_2} \cdots \gamma_k^{m_k} \gamma_1^{-n_1} \gamma_2^{-n_2} \cdots \gamma_k^{-n_k}  $$
where $m_1 + \cdots + m_k + n_1 + \cdots + n_k = m-1$. Notice that if $m_k$ and $n_1$ are both nonzero, then no occurrence of~$\gamma_i$ is immediately followed by~$\gamma_i^{-1}$; so we have $\gamma_{i+1}^* = \gamma_i^*$ whenever $j_{i+1} = j_i$, as desired. Therefore, by permuting $\gamma_1,\ldots,\gamma_k$, we may assume $m_i = n_i = 0$ for all $i > 1$. Also, we may assume $m_1$ and~$n_1$ are both nonzero, for otherwise we have $\gamma_i^* = \gamma_j^*$ for all~$i$ and~$j$. Then, since $\gamma_1 \gamma_1^{-1} = \gamma_2 \gamma_2^{-1}$, we have
	$$  \gamma_1^* \gamma_2^* \cdots \gamma_{m-1}^*
	= \gamma_1^{m_1} \gamma_1^{-n_1} 
	= \gamma_1^{m_1-1}\gamma_2 \gamma_1^{-(n_1-1)} \gamma_2^{-1} . $$
Assuming, without loss of generality, that $n_1 \ge m_1$, so $n_1 \ge \lceil (m-1)/2 \rceil \ge 2$,  this new representation of the same product satisfies the condition that $\gamma_i$ is never immediately followed by~$\gamma_i^{-1}$. This completes the proof.
\end{proof}

\subsection{Facts from group theory}

\begin{lem} \label{G'Divides}
If $|G_k'|$ is square-free, then $|G_k' / G_{k-1}'|$ is a divisor of both $|\quot{G_{k-1}}|$ and $|\quot{G_k}/\quot{G_{k-1}}|$.
\end{lem}

\begin{proof} 
We may assume $k = \ell$, so $G = G_k$. Let $p$ be a prime factor of $|G' / G_{k-1}'|$, let $P$ be the Sylow $p$-subgroup of~$G$, and let $\varphi \colon G \to P$ be the natural projection. Since $|G'|$ is square-free, it suffices to show that $|\quot{G_{k-1}}|$ and $|\quot{G_k}/\quot{G_{k-1}}|$ are divisible by~$p$.

We may assume $|G'| = p$ and $G_{k-1}' = \{e\}$, by modding out the unique subgroup of index~$p$ in~$G'$. Therefore $\varphi (G_{k-1})$ is abelian, so it is a proper subgroup of~$P$. Since $G' = P' \subset \Phi(P)$, this implies $\varphi (G_{k-1}) G'$ is a proper subgroup of~$P$, so its index is divisible by~$p$. Hence $|\quot{G} / \quot{G_{k-1}}|$ is divisible by~$p$.

There must be some $t \in S_{k-1}$, such that $[\sigma_k,t]$ is nontrivial. Hence $\varphi(t) \notin Z(G) \supset G'$, so $p$ is a divisor of $|\quot{\varphi(t)}|$, which is a divisor of~$|\quot{G_{k-1}}|$.
\end{proof}

The following fact is well known and elementary, but we do not know of a reference in the literature. It relies on our assumption that $G'$ is cyclic.

\begin{lem} \label{CommsGen} 
We have $\langle\, [s,t] \mid s,t \in S \,\rangle = G'$ if $N \subset Z(G)$.
\end{lem}

\begin{proof}
Let $H = \langle\, [s,t] \mid s,t \in S \,\rangle$. Then $H$ is a normal subgroup of~$G$, because every subgroup of a cyclic, normal subgroup is normal. In $G/H$, every element of~$S$ commutes with all of the other elements of~$S$ (and with all of~$N$), so $G/H$ is abelian. Hence $G' \subset H$.
\end{proof}

\section{Base case of the inductive construction} \label{BaseCaseSect}

Recall that the condition $\alpha_k^\epsilon$  is defined in \cref{AssumpsSect}.

\begin{prop}[(cf.\ {\cite[Case~6.2]{KeatingWitte}})] \label{KW6-2}
Assume $|N|$ is square-free (and $\ell \ge 2$).
Then $\alpha_2^2$ is true unless $|G_2'| = m_2 =|\quot{\sigma_1}| = |\quot{\sigma_2}| = 3$.
\end{prop}

\begin{proof}
For convenience, let
	$$ \text{$a = \sigma_1$, \ $b = \sigma_2$, \ and \ $m = m_2$,} $$
and define $r$ by 
	$$ \text{$\quot{b}^{m} = \quot{a}^r $ \ and \ $0 < r \le |\quot{a}|$.} $$
We may assume:
	\noprelistbreak
	\begin{itemize}
	\item $\ell = 2$, so $S = S_2 = \{a,b\}$ and $G = G_2$.
	\item $(G')^2$ is nontrivial. (Otherwise, the condition about generating $(G')^2$ is automatically true, so it suffices to show $\Voltages{2} \neq \emptyset$, which is easy.)
	\item Either $|\quot{a}|$ is even, or  $m$ is odd (by interchanging $\sigma_1$ and~$\sigma_2$ if necessary).
	\item $|\quot{a}| \neq 3$ (by interchanging $\sigma_1$ and~$\sigma_2$ if necessary: if $|\quot{\sigma_1}| = |\quot{\sigma_2}| = 3$, then $m = 3$ and, from \cref{G'Divides}, we also have $|G'| = 3$, which means we are in a case in which the statement of the \namecref{KW6-2} does not make any claim).
	\item $r \ge |\quot{a}|/2$ (by replacing $a$ with its inverse if necessary).
	\end{itemize}

Note that $|G'|$ is a divisor of both $|\quot{a}|$ and~$m$ \csee{G'Divides}. Since $(G')^2$ is nontrivial, this implies that $|\quot{a}|$ and~$m$ both have at least one odd prime divisor.

\setcounter{case}{0}

\begin{case} \label{BaseCase(m=3)}
Assume $m = 3$.
\end{case} 
Since $|G'|$ is a divisor of $m$, we must have $|G'| = 3$, so $|\quot{a}|$ must be divisible by~$3$. 	Then, since $|\quot{a}| \neq 3$, we must have $|\quot{a}| \ge 6$. Furthermore, by applying \cref{G'Divides} with $a$ and~$b$ interchanged, we see that $|\quot{G}/\langle \quot{b} \rangle|$ is also divisible by $|G'| = 3$, which means that $r$~is divisible by~$3$. 

We claim that it suffices to find two elements $\gamma_1,\gamma_2 \in \Voltages{2}$, such that $\gamma_1 \neq \gamma_2$ and $\gamma_1 \in \gamma_2 G'$. To see this, note that, for any $x \in N$, there is some $i \in \{1,2\}$, such that $\langle \gamma_i x \rangle$ has nontrivial projection to~$G'$ (with respect to the unique direct-product decomposition $N = G' \times (G')^\perp$). Since $|G'|$ is prime, this implies that the projection is all of~$G'$, so \cref{ContainsProj} tells us that $\langle \gamma_i x \rangle$ contains~$G'$. This establishes $\alpha_2^1$, which is equivalent to~$\alpha_2^2$ \fullcsee{Alpha<>Alpha}{odd}.
This completes the proof of the claim.

Assume, for the moment, that $r = 3$. Then, since $r \ge |\quot{a}|/2$ and $|\quot{a}| \ge 6$, we must have $|\quot{a}| = 6$. Here are two hamiltonian cycles in $\Cay( \quot{G}; \quot{a}, \quot{b} )$ that cover~$S^{\pm1}$:
	$$ ( b^{-1}, a^{-2}, b^{-4}, a^{-2}, b^{-1}, a^3, b^2, a, b^{-2} ) $$
and
	$$ ( b^{-1}, a^{-2}, b^{-1}, a, b^{-1}, a^{-1}, b^{-2}, a^{-1}, b^{-1}, a^2, b^2, a, b^{-2} ) $$
\csee{Case6x3Fig}.
Straightforward calculations show that their voltages are $b^{-6} [a,b]$ and $b^{-6} [a,b]^2$, respectively. So we may let $\gamma_1 = b^{-6} [a,b]$ and $\gamma_2 = b^{-6} [a,b]^2$.

\begin{figure}[t]
\begin{center}
\includegraphics{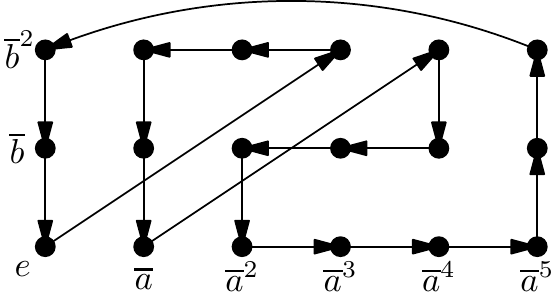}
\hfil
\includegraphics{PDF/6x3-1.pdf}
\caption{Two hamiltonian cycles in $\Cay \bigl( \quot{G}; \{ \quot{a}, \quot{b} \} \bigr)$
when $m = r = 3$.}
\label{Case6x3Fig}
\end{center}
\end{figure}

We may now assume $r \ge 6$ (since $r$~is divisible by~$3$). Let
	$$ I = \begin{cases}
	\{0,1\} & \text{if $r \neq |\quot{a}|$} , \\
	\{1,2\} & \text{if $r = |\quot{a}|$} . \\
	\end{cases} $$
Then, for $i \in I$, we have $0 \le i \le |\quot{a}| - 4$, and $4 \le r - i \le |\quot{a}|-1$. So the walk
	$$ C_i = \bigl( a^i, b^{-1}, a^{-(|\quot{a}| - r + i -1)}, 
	b^{-1}, a^{|\quot{a}| - 4}, b^{-1},
	a^{-(|\quot{a}|-i-4)}, b^{-1}, 
	a^{r-i-3}, b^{-2}, a, b^2, a,
	b^{-2} \bigr) $$
is as pictured in \cref{Case4-3m=3Fig}.
It is a hamiltonian cycle in $\Cay(\quot{G}; \quot{a},\quot{b})$ that covers $S^{\pm1}$, and its voltage is of the form $[a,b]^{-2i} h_2$, where $h_2$~is independent of~$i$. Thus, we may let 
	$$\{\gamma_1,\gamma_2\} = \{\, [a,b]^{-2i} h_2 \mid i \in I \,\} . $$
\begin{figure}[t]
\begin{center}
\includegraphics{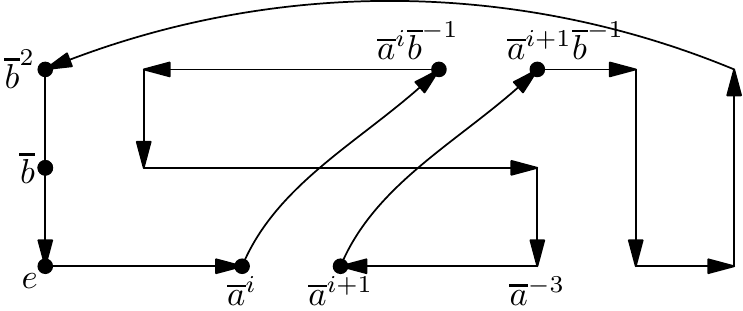}
\caption{A hamiltonian cycle in $\Cay \bigl( \quot{G}; \{ \quot{a}, \quot{b} \} \bigr)$
when $m = 3$ and $r \ge 6$.}
\label{Case4-3m=3Fig}
\end{center}
\end{figure}
%
%
%
%
%
%
%
%
%
%
%
%
%
%
%

\begin{case} \label{Case4-3MOdd}
Assume $m \neq 3$.
\end{case}
(Cf.\ \cite[Case~4.3]{KeatingWitte}.)
Since $m$ and $|\quot{a}|$ both have at least one odd prime divisor, we must have $m \ge 5$ and $|\quot{a}| \geq 5$. 
Let 
	$$ X = 
	\begin{cases}
	\bigl( b^{-(m-2)}, a, b^{m-3}, a^{|\quot{a}| - 3}, b^{-1}, ( a^{-(|\quot{a}| - 4)}, b^{-1}, a^{|\quot{a}| - 4}, b^{-1} )^{(m-3)/2} \bigr) 
		& \text{if $|\quot{a}|$ is odd,} \\
	\hfil \bigl( b^{-1}, 
			(b^{-(m-3)}, a, b^{m-3}, a)^{(|\quot{a}|/2) - 1},
			b^{-(m-2)}	
			\bigr) 
		& \text{if $|\quot{a}|$ is even}
		.\end{cases}
		$$
For each~$i$ with $1 \le i \le \lfloor(|\quot{a}| - 1)/2\rfloor$, we have
$1 \le i \le \min \bigl(r-1, |\quot{a}|-3 \bigr)$ (since $r \ge |\quot{a}|/2$ and $|\quot{a}| \ge 5$),
so we may let
	$$ C_i = \bigl( a^i, b^{-1}, a^{-(|\quot{a}|+i-r-1)}, X, a^{-(|\quot{a}|-i-2)}, b^{-1}, a^{r-i-1}, b^{-(m-1)} \bigr) $$
\csee{Case4-3OddFig,Case4-3EvenFig}.
Then $C_i$ is a hamiltonian cycle in $\Cay(\quot{G}; \quot{a},\quot{b})$.

Note that both possibilities for~$X$ contain oriented edges labelled $a$, $b$, and~$b^{-1}$. Furthermore, since $|\quot{a}| - i -2 \ge 1$, we see that $C_i$ also contains at least one oriented edge labelled~$a^{-1}$. Therefore $C_i$ covers $\{a, b, a^{-1}, b^{-1}\} = S^{\pm1}$. 

The voltage $\voltage C_i$ of~$C_i$ is of the form $[a,b]^{-2i} h_2$, where $h_2$~is independent of~$i$.
Since $|\quot{a}| \ge |G'|$ \csee{G'Divides} and $\langle [a,b] \rangle = G'$ \csee{CommsGen},  \cref{GCD(consecutive)} (combined with \cref{ContainsProj}) tells us that for any $x \in N$, we may choose~$i$ so that $\langle (\voltage C_i)x \rangle$ contains $(G')^2$.
\end{proof}

\begin{figure}[t]
\begin{center}
\includegraphics{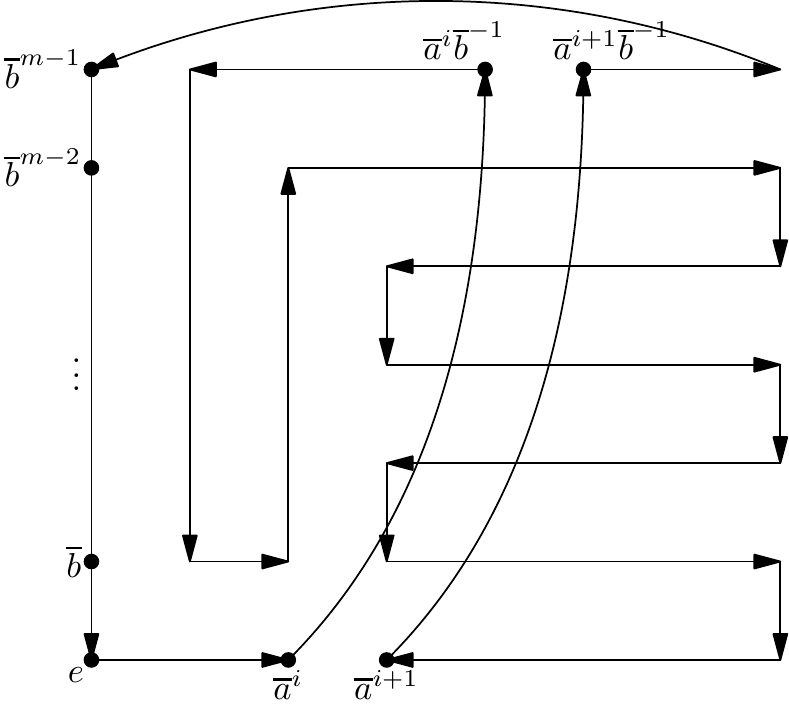}
\caption{A hamiltonian cycle~$C_i$ in $\Cay \bigl( \quot{G}; \{ \quot{a}, \quot{b} \} \bigr)$ when $m = |\quot{G}/\langle \quot{a} \rangle|$ is odd.}
\label{Case4-3OddFig}
\end{center}
\end{figure}
%
%
%
%
%
%
%
%
%
%
%
%
%
%
%

\begin{figure}[t]
\begin{center}
\includegraphics{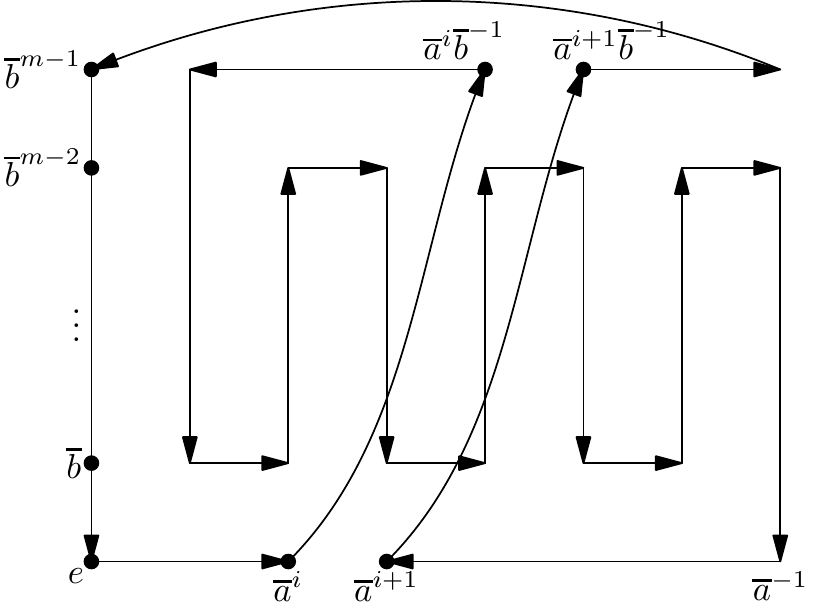}
\caption{A hamiltonian cycle $C_i$ in $\Cay \bigl( \quot{G}; \{ \quot{a}, \quot{b} \} \bigr)$
when $|\quot{a}|$ is even.}
\label{Case4-3EvenFig}
\end{center}
\end{figure}
%
%
%
%
%
%
%
%
%
%
%
%
%
%
%

\section{The main induction step} \label{InductionStepSect}

\begin{figure}[ht]
\begin{center}
\includegraphics{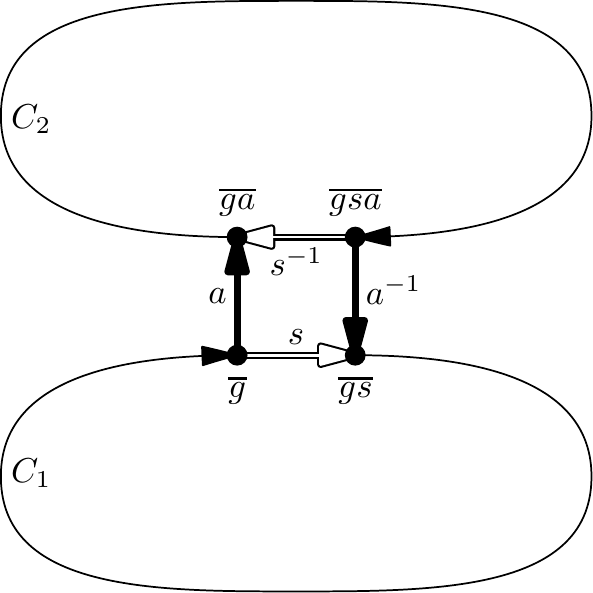}
\end{center}
\caption{$C_1$ and~$C_2$ are merged into a single cycle by replacing the two white edges labelled $s$ and~$s^{-1}$ with the two black edges labelled $a$ and~$a^{-1}$.}
\label{ConnSumFig}
\end{figure}

The induction step of our proof uses the following well-known gluing technique that is illustrated in \cref{ConnSumFig}. 
 
 \begin{defn} \label{ConnectedSumDef}
 Let 
 	\begin{itemize}
	\item $C_1$ and~$C_2$ be two disjoint oriented cycles in $\Cay(\quot{G}; \quot{S})$,		\item $g \in G$,
	and
	\item $a, s \in S$.
	\end{itemize}
If 
	\begin{itemize}
	\item $C_1$ contains the oriented edge $\quot{g}(s)$,
	and 
	\item $C_2$ contains the oriented edge $\quot{gsa}(s^{-1})$,
	\end{itemize}
then we use $C_1 \connsum_s^a C_2$ to denote the oriented cycle obtained from $C_1 \cup C_2$ as in \cref{ConnSumFig}, by
	\begin{itemize}
	\item removing the oriented edges $\quot{g}(s)$ and $\quot{gsa}(s^{-1})$, and
	\item inserting the oriented edges $\quot{g}(a)$ and $\quot{gsa}(a^{-1})$.
	\end{itemize}
This may be called the \emph{connected sum} of~$C_1$ and~$C_2$.
 \end{defn}

%
%
%
%
%
%
%
%
%
%

 \begin{lem} \label{EndptOfConnSum}
If $C_1$, $C_2$, $g$, $s$, and~$a$ are as in \cref{ConnectedSumDef}, 
and $N \subset Z(G)$, then
	$$ \voltage (C_1 \connsum_s^a C_2) = (\voltage C_1) (\voltage C_2) [a,s] .$$
\end{lem}

\begin{proof}
Write $C_1 = \quot{gs}(s_i)_{i = 1}^m$ and $C_2 = \quot{ga}(t_j)_{j=1}^n$. Then
	$$  C_1 \connsum_s^a C_2 = \quot{gsa} \bigl( a^{-1}, (s_i)_{i = 1}^{m-1}, a, (t_j)_{j=1}^{n-1} \bigr) ,$$
so 
	\begin{align*}
	 \voltage (C_1 \connsum_s^a C_2)
	 &= a^{-1} \, \left( \prod_{i = 1}^{m-1} s_i \right) \, a \, \left( \prod_{j = 1}^{n-1} t_j \right) 
	 \\&=a^{-1} \,  \left( \prod_{i = 1}^{m} s_i \right) \! s_m^{-1} \ a \, \left( \prod_{j = 1}^{n} t_j \right) \! t_n^{-1} 
	 \\&= a^{-1} \, \left(\voltage C_1 \right) s^{-1} \, a \, \left( \voltage C_2 \right) s 
	 \\&= \left(\voltage C_1 \right)  \left( \voltage C_2 \right) \, a^{-1} s^{-1}  a  s 
	 && (\voltage C_i \in N \subset Z(G))
	 \\&= \left(\voltage C_1 \right)  \left( \voltage C_2 \right) \, [a,s]
	. &&\qedhere \end{align*}
\end{proof}

\begin{cor} \label{CanGlueCycles}
Assume 
	\begin{itemize}
	
	\item $2 \le k < \ell$, and\/ {\upshape(}to eliminate some subscripts\/{\upshape)} $m = m_{k+1}$ and $a = \sigma_{k+1}$,
	
	\item $\pi_1,\pi_2,\ldots,\pi_m$ are elements of~$\Voltages{k}$,
	
	\item $s_1,s_2,\ldots,s_{m-1}$ are elements of $S_k$, and, for each~$i$, a choice $s_i^* \in \{s_i^{\pm1} \}$ has been made in such a way that if $s_{i+1} = s_i$, then $s_{i+1}^* = s_i^*$,
	and

	\item $N \subset Z(G)$.
	\end{itemize}
Then there is a hamiltonian cycle in $\Cay( \quot{G_{k+1}}; \quot{S_{k+1}})$ that  covers~$S_{k+1}^{\pm1}$, and whose voltage is
	$$ \left( \prod_{i=1}^{m} \pi_i \right) \left( \prod_{i=1}^{m-1} [a, s_i^*] \right).$$
\end{cor}
 
 \begin{proof}
For each~$i$, let $C_i$ be an oriented hamiltonian cycle in $\Cay( \quot{G_k}; \quot{S_k} )$ that  covers $S_k^{\pm1}$, and has voltage~$\pi_i$. We inductively construct sequences $(g_i)_{i=1}^m$ and $(x_i)_{i=1}^m$ of elements of~$G_k$, as follows.

Let $g_1 = e$. Since $C_1$ covers $S_k^{\pm1}$, we know there is some $x_1 \in G_k$, such that $ag_1C_1$ contains the oriented edge $\quot{ax_1}(s_1^*)$.

Now, suppose $g_1, x_1, g_2, x_2, \ldots, g_i, x_i \in G_k$ are given, such that the connected sum
	$$a g_1 C_1 \connsum_{s_1^*}^a a^2 g_2 C_2  \connsum_{s_2^*}^a  \cdots  \connsum_{s_{i-1}^*}^a a^{i} g_{i} C_{i} $$
exists, and contains the oriented edge $\quot{a^{i} x_{i}}(s_{i}^*)$. Since $C_{i+1}$ covers $S_k^{\pm1}$, we know that $C_{i+1}$ contains an oriented edge labelled $(s_{i}^*)^{-1}$, and a different oriented edge that is labelled~$s_{i+1}^*$. Therefore, there exist $g_{i+1},x_{i+1} \in G_k$, such that 
	$$ \text{$a^{i+1} g_{i+1} C_{i+1}$ contains the oriented edges $\quot{a^{i+1} x_{i}s_{i}^*}\bigl( (s_{i}^*)^{-1} \bigr)$ and  $\quot{a^{i+1} x_{i+1}}( s_{i+1}^*)$.} $$
The first of these edges is removed when we form the connected sum 
		$$\bigl( a g_1 C_1 \connsum_{s_1^*}^a a^2 g_2 C_2  \connsum_{s_2^*}^a  \cdots  \connsum_{s_{i-1}^*}^a a^{i} g_{i} C_{i} \bigr) 
		\  \connsum_{s_{i}^*}^a \ 
		a^{i+1} g_{i+1} C_{i+1} ,$$
 but the second edge remains, and will be used to form the next connected sum (unless $i +1 = m$).
 
 Since each $C_i$ is a hamiltonian cycle in $\Cay( \quot{G_k}; \quot{S_k})$, the resulting connected sum 
 	$$a g_1 C_1 \connsum_{s_1^*}^a a^2 g_2 C_2  \connsum_{s_2^*}^a  \cdots  \connsum_{s_{m-1}^*}^a a^{m} g_{m} C_{m} $$
passes through all of the vertices in $\quot{a G_k} \cup \quot{a^2 G_k} \cup \cdots \cup \quot{a^m G_k}$. That is, it passes through every element of~$\quot{G_{k+1}}$, so it is a hamiltonian cycle.
Its voltage is calculated by repeated application of \cref{EndptOfConnSum}.

To complete the proof, we verify that the hamiltonian cycle covers $S_{k+1}^{\pm1}$.
Since each $C_i$ covers $S_k^{\pm1}$, the disjoint union 
	$$a g_1 C_1 \cup a^2 g_2 C_2  \cup  \cdots  \cup a^{m} g_{m} C_{m} $$
contains (at least) $m$ disjoint pairs of edges labelled $s$ and~$s^{-1}$, for each $s \in S_k$. Each invocation of the connected sum removes only one such pair, and the operation is performed only $m-1$~times, so at least one of the $m$ pairs must remain, for each $s \in S_k$. Therefore, the hamiltonian cycle covers $S_k^{\pm1}$. Also, the cycle certainly covers~$a^{\pm1}$, since each invocation of the connected sum inserts a pair of edges labelled $a$ and~$a^{-1}$. Hence, the hamiltonian cycle covers $S_k^{\pm1} \cup \{a^{\pm1}\} = S_{k+1}^{\pm1}$.
 \end{proof}

We can now prove the main result of this section. (Recall that the condition $\alpha_k^\epsilon$  is defined in \cref{AssumpsSect}.)

\begin{prop} \label{Alpha->Alpha}
Assume $|N|$ is square-free and $|G_{k+1}'/G_k'|$ is odd. Then
	\begin{enumerate}
	\item  \label{Alpha->Alpha-1}
	$\alpha_k^1 \Rightarrow \alpha_{k+1}^1$,
	and
	\item  \label{Alpha->Alpha-2+}
	$\alpha_k^{2+} \Rightarrow \alpha_{k+1}^{2+}$
	if $|[s,t]|$ is even for all $s,t \in S_{k+1}$ with $s \neq t$.
	\end{enumerate}
\end{prop}

\begin{proof}
For convenience, let $m = m_{k+1}$ and $a = \sigma_{k+1}$.
Choose an oriented hamiltonian cycle $C$ in $\Cay(\quot{G_k}; \quot{S_k})$ that covers $S_k^{\pm1}$, and has its endpoint in $h_k (G_k')^\epsilon$. There is no harm in assuming that the endpoint is precisely~$h_k$. Let
	$$ h_{k+1} = (h_k)^m [a,\sigma_1]^{m-1}.$$

Since $\{\, [s,t] \mid s,t \in S_{k+1} \,\}$ generates~$G_{k+1}'$ \csee{CommsGen}, we know that $\{\, [a,s] \mid s \in S_k \,\}$ generates $G_{k+1}'/G_k'$. Therefore, given any $x \in N$, \cref{ContainsN2}  (combined with \cref{ContainsProj}) tells us we may choose a sequence $(s_i)_{i=1}^{m-1}$ of elements of $S_k$, and $s_i^* \in \{s_i^{\pm1}\}$ for each~$i$, such that  $s_{i+1}^* = s_i^*$ whenever $s_{i+1} = s_i$, and
	\begin{align} \label{WithStarGensG'}
	\text{$\left\langle x \, (h_k)^m  \prod_{i=1}^{m-1} [a,s_i^*] ,
	\ (G_k')^\epsilon \right\rangle$ contains $(G_{k+1}')^\epsilon$}
	 .\end{align}
From $(\alpha_k^\epsilon)$, we know there exists $\pi \in \Voltages{k} \cap h_k \,  (G_k')^\epsilon$, such that, if we let
	$$ \gamma = \pi \, (h_k)^{m-1} \,  \prod_{i=1}^{m-1} [a, s_i^*] , $$
then $\langle x \gamma \rangle$ contains~$(G_k')^\epsilon$. Since $\pi \equiv h_k \pmod{(G_k')^\epsilon}$, combining this with \pref{WithStarGensG'} shows that $\langle x \gamma \rangle$ contains~$(G_{k+1}')^\epsilon$. Also, since we are assuming $|[a, s_i^*]|$ is even if $\epsilon = 2$, we have $[a, s_i^*] \equiv [a, \sigma_1] \pmod{(G_{k+1}')^\epsilon}$ for all~$i$, so
	$$\gamma \in (h_k)^m [a, \sigma_1]^{m-1} (G_{k+1}')^\epsilon = h_{k+1} (G_{k+1}')^\epsilon .$$

Furthermore, \cref{CanGlueCycles} tells us that there is a hamiltonian cycle in $\Cay( \quot{G}; \quot{S} )$ whose voltage is~$\gamma$, and this hamiltonian cycle covers $S^{\pm1}$.
This establishes $\alpha_{k+1}^\epsilon$.

Now, if $\epsilon = 2$, then our assumptions imply that $|h_k|$ and $|[a,\sigma_1]|$ are both even. Since $m$ and~$m-1$ are of opposite parity, this implies that $|h_{k+1}|$ is even, so $\langle h_k, (G'_{k+1})^2 \rangle$ contains~$G'_{k+1}$. This establishes $\alpha_{k+1}^{2+}$.
\end{proof}

\section{Combining the base case with the induction step} \label{CombiningSect}

Recall that the condition $\alpha_k^\epsilon$  is defined in \cref{AssumpsSect}.

\begin{cor} \label{SpiralG'Odd}
Assume $|N|$ is square-free and $\ell \ge 3$. If $|G'|$ is odd, then $\alpha_\ell^1$ is true unless $|G'| = |\quot{s}| = 3$ for all $s \in S$.
\end{cor}

\begin{proof}
Assume it is not the case that $|G'| = |\quot{s}| = 3$ for all $s \in S$.
Then we may assume  (by permuting the elements of~$S$) that either $|G_2'| \neq 3$ or $|\sigma_1| \neq 3$. Therefore \cref{KW6-2} tells us that $\alpha_2^2$ is true.
Also, since $|G'|$ is odd, we have $\alpha_2^2 \Leftrightarrow \alpha_2^1$ \fullcsee{Alpha<>Alpha}{odd}, so $\alpha_2^1$ is true. Then repeated application of \fullcref{Alpha->Alpha}{1} establishes~$\alpha_\ell^1$.
\end{proof}

\begin{prop} \label{SpiralG'Even}
Assume $|N|$ is square-free and $\ell \ge 3$. If $|G'|$ is even, then:
	\begin{enumerate}
	\item $\alpha_\ell^1$ is true if there exist $s,t \in S$, such that $|[s,t]|$ is odd and $s \neq t$.
	\item $\alpha_\ell^{2+}$ is true if $|[s,t]|$ is even for all $s,t \in S$ with $s \neq t$.
	\end{enumerate}
\end{prop}

\begin{proof}
Since $|G'|$ is even, we may assume (by permuting the elements of~$S$) that $|[\sigma_3, \sigma_1]|$ is even. It suffices to prove $\alpha_3^1$ or $\alpha_3^{2+}$ (as appropriate), for then repeated application of \cref{Alpha->Alpha} establishes the desired conclusion. Thus, we may assume $\ell = 3$, so $G_3 = G$. Let $m = m_3$ and $a = \sigma_3 = \sigma_\ell$. 

By permuting the elements of~$S$, we may assume that either:
	\begin{enumerate} \itemindent=1in
	\item[odd case:] $|[\sigma_3, \sigma_2]|$ is odd,
	or
	\item[even case:] $|[s,t]|$ is even for all $s,t \in S$ with $s \neq t$.
	\end{enumerate}
Furthermore, in the even case, we may assume that either:
	\begin{enumerate} \itemindent=1in
	\item[even  subcase:] $|[\quot{\sigma_1}, \quot{\sigma_2}]|$ has even index in~$\quot{G}$,
	or
	\item[odd subcase:] $\langle \quot{s}, \quot{t} \rangle$ has odd index in~$\quot{G}$, for all  $s,t \in S$, such that $s \neq t$.
	\end{enumerate}

Since $|[\sigma_3, \sigma_1]|$ is even, we know $|\quot{\sigma_1}|$ is even \csee{G'Divides}, so $|\quot{\sigma_1}| \neq 3$. Therefore \cref{KW6-2} tells us that $\alpha_2^2$ is true.

We now use a slight modification of the proof of \cref{Alpha->Alpha}.
Choose an oriented hamiltonian cycle $C$ in $\Cay(\quot{G_2}; \quot{S_2})$ that covers $S_2^{\pm1}$, and has its endpoint in $h_2 (G_2')^2$. There is no harm in assuming that the endpoint is precisely~$h_2$. 

\cref{TwoGen} provides a hamiltonian cycle $C'$ in $\Cay(\quot{G_2}; \quot{S_2})$, such that $|\voltage C'|$ is even. Let
	$$ h' =
	\begin{cases}
	\voltage C' & \text{in the odd subcase of the even case}, \\
	\hfil h_2 & \text{in all other cases}.
	\end{cases}
	$$
Let $h_3 = (h_2)^{m-1} h' [a,\sigma_1]^{m-1}$.

Since $\{\, [s,t] \mid s,t \in S \,\}$ generates~$G'$ \csee{CommsGen}, we know that $\{\, [a,s] \mid s \in S_2 \,\}$ generates $G'/G_2'$. 
Therefore, given any $x \in N$, \cref{ContainsN2}  (combined with \cref{ContainsProj}) tells us we may choose a sequence $(s_i)_{i=1}^{m-1}$ of elements of $S_2$, and $s_i^* \in \{s_i^{\pm1}\}$ for each~$i$, such that $s_{i+1}^* = s_i^*$ whenever $s_{i+1} = s_i$, and
	\begin{align} \label{3WithSquareGensG'}
	\text{$\left\langle x (h_2)^{m-1}h' \prod_{i=1}^{m-1} [a,s_i^*] , 
	\ (G_2')^2 \right\rangle$ contains $(G')^2$}
	. \end{align}
Furthermore, in the odd case, the choices can be made so that \pref{3WithSquareGensG'} holds with~$G'$ in the place of~$(G')^2$.

From  $\alpha_2^2$, we know there exists $\pi \in \Voltages{2} \cap h_2 \,  (G_2')^2$, such that, if we let
	$$ \gamma = \pi \, (h_2)^{m-2}h' \,  \prod_{i=1}^{m-1} [a,s_i^*] , $$
then 
	\begin{align} \label{3SquareInGamma}
	\text{$\langle x\gamma \rangle$ contains~$(G')^2$.}
	\end{align}
It is clear from the definitions that $\gamma \in h_3 G_3'$. Furthermore, we have $\gamma \in h_3 (G_3')^2$ in the even case.

\cref{CanGlueCycles} tells us that there is a hamiltonian cycle in $\Cay( \quot{G}; \quot{S} )$ whose voltage is~$\gamma$, and this hamiltonian cycle covers $S^{\pm1}$. We now consider various cases individually.

\setcounter{case}{0}

\begin{case}
The odd case.
\end{case}
Recall that, in this case, \pref{3WithSquareGensG'} holds with $G'$ in the place of~$(G')^2$. 
Since $\pi \equiv h_2 \pmod{(G_2')^2}$, combining this with \pref{3SquareInGamma} shows that $\langle x \gamma \rangle$ contains all of~$G'$. This establishes $\alpha_3^1$.

\begin{case}
The even subcase of the even case.
\end{case}
In this subcase, we know $m$ is even, $h' = h_2$, and $|[a,\sigma_1]|$ is even. Since $h_3 = (h_2)^m [a,\sigma_1]^{m-1}$, we see that $|h_3|$ is even, so $\langle h_3, (G_3')^2 \rangle$ contains~$G_3'$. This establishes $\alpha_3^{2+}$.

\begin{case}
The odd subcase of the even case.
\end{case}
In this subcase, we know $m-1$ is even, and $h' = \voltage C'$ has even order. Therefore $|h_3|$ is even, so $\langle h_3, (G_3')^2 \rangle$ contains~$G_3'$. This establishes $\alpha_3^{2+}$.
\end{proof}

\end{document}